\documentclass[12pt]{article}
\newtheorem{df}{Definition}[section]
\newtheorem{definition}[df]{Definition}
\newtheorem{theorem}[df]{Theorem}

\newtheorem{corollary}[df]{Corollary}

\begin{document}
\title{On a Generalization of Property B}

\author{Robert Cowen\\Mathematics Department\\Queens College, CUNY\\Flushing, 
NY 11367}

\date{\today}
\maketitle

\begin{abstract}
A set-theoretic property called Property S  is introduced as a generalization of the well-known Property B . Property S is named for  Schrijver who first used it in a paper \cite{Sch} in connection with the Boolean prime ideal theorem. It was independently introduced by Kolany \cite{Kol} to give a uniform treatment for a variety of satisfiability problems who then used a generalized resolution method to determine satisfiability. Here we further investigate Property S and the resolution method.

\end{abstract}
\section{Introduction}
A collection $\cal{E}$ of subsets of a set $V$ is said to have {\bf Property B} if there exits a partition, \{$X$,$V-X$\}, of $V$,  such that both $X$, $V-X$ intersect every $E$ in $\cal{E}$. A set which intersects every set in a family of sets is often called a {\bf transversal} for that family; thus we can simply say: $\cal{E}$ has Property B if there exists a subset $X$ such that both $X$, $V-X$ are transversals for $\cal{E}$. Clearly $V-X$ is a transversal for $\cal{E}$ if and only if $X$ does not contain any $E$ in $\cal{E}$. Hence, $\cal{E}$ has Property B if and only if there is a transversal, $X$, for $\cal{E}$ which does not contain any $E$ in $\cal{E}$. We shall also say that the hypergraph $<V, \cal{E}>$ has Property B if $\cal{E}$ has the property.

The property is named for Felix Bernstein who proved, in 1908, that a countable system of infinite sets has Property B. Since then it has been studied extensively (see, for example, \cite{Erd1},\cite{Erd2}, \cite{Erd3}, \cite{EH}, \cite{Mil}). In graph theory Property B has been linked by Woodall \cite{Woo} and Stein \cite{Ste} with the 4-color problem. For example Woodall\cite {Woo}  showed that the 4-color conjecture was equivalent to the odd circuits of a planar graph, when regarded as sets of edges, has Propery B. 

Property B is an NP-Complete property when restricted to finite families of finite sets; in fact, it remains NP-Complete even when all sets in the family have at most three elements (see \cite{Cow93}). In the next section we condsider a generalization to two families of subsets and we name the generalization ``Property S" for reasons explained below. Property B is then the special case when both families are the same.

\section{Property S}
Two families $\cal{E}$, $\cal{F}$ of subsets of $V$ will be said to have {\bf Property S} if there exists a partition, \{ $X$, $V-X$ \}, of $V$,  such that $X$ is a transversal for $\cal{E}$ and $V-X$ is a transversal for $\cal{F}$; the partition will be referred to as an {\bf S-partition}. Equivalently, $\cal{E}$, $\cal{F}$ have Propery S if there is a transversal $X$ for $\cal{E}$, which does not contain any of the subsets of $\cal{F}$. 

If $\cal{E}$, $\cal{F}$ has Property S, then $\cal{E} \cap \cal{F}$ has Property B, since every set in $\cal{E} \cap \cal{F}$ must then intersect both cells of the partition. If $\cal{E} =\cal{F}$,  $\cal{E}$, $\cal{F}$ has Property S if and only if $\cal{E}$ has Property B. 

From now on we shall often write $<V,\cal{E}, \cal{F}>$ where  $\cal{E}$, $\cal{F}$ are collections of subsets of $V$ and ask whether this ``bihypergraph'' has Property S. Instead of explicitly mentioned the partition of $V$, Schrijver \cite{Sch} refers to disjoint transversals for the families $\cal{E}$, $\cal{F}$, but this is clearly equivalent. He proved the following ``compactness'' theorem for property S and showed that it is equivalent in ZF-set theory to BPI, the prime ideal theorem for Boolean algebras. 

\begin{theorem} Let  $\cal{E}$, $\cal{F}$ be families of finite subsets of $V$. Then $<V,\cal{E}, \cal{F}>$ has Property S if $<V,\mathcal{E}_0, \mathcal{F}_0>$ has Property S for every finite $\mathcal{E}_0\subset\cal{E}$, $\mathcal{F}_0\subset \cal{F}$.
\end{theorem}

Restrictions of Schrijver's Theorem are also equivalent to BPI; for example,  the corresponding compactness result for Property B is equivalent to BPI even if all sets in the family have at most three elements (see \cite{Cow93}). 

As mentioned in the Introduction, Bernstein showed that a countable system of infinite sets has Property B; this can be generalized as follows.
\begin{theorem}
Let $\mathcal{E}$, $\mathcal{F}$ be countable collections of infinite sets and let $V=(\cup \mathcal{E}) \cup  (\cup \mathcal{F})$. Then $<V,\cal{E}, \cal{F}>$ has Property S.
\end{theorem}
{\bf Proof}. Suppose $\mathcal{E}=\{A_n\}^\infty_{n=1}$, $\mathcal{F}=\{B_n\}^\infty_{n=1}$. For $n\ge 1$, choose elements $x_n, y_n$ as follows:\\
 $x_1 \in A_1$; $y_1 \in B_1$, $y_1 \neq x_1$; 
$x_{n+1} \in A_{n+1}$, $x_{n+1} \neq y_i$, $1 \leq i \leq n$; 
$y_{n+1} \in  B_{n+1}$, $y_{n+1} \neq x_i$, $1 \leq i \leq n+1$.  Let $X={\cup }^\infty_{n=1} \{x_n\}$, $Y={\cup }^\infty_{n=1} \{y_n\}$. Then $X$,$Y$ are disjoint transversals for $\mathcal{E}$, $\mathcal{F}$, respectively and the Theorem easily follows.\bigskip

Much is known about Property B in the finite case, as well (see \cite{Erd1}, \cite{Erd2},\cite{Erd3},\cite{Mil},\cite{Woo}). For example, it is easy to show that any family of subsets of a $(2n+1)$ element set $V$ fails to have Property B if it contains all of the $(n+1)$ element subsets of $V$. We can generalize this as follows.
\begin{theorem}
Suppose $V$ is a set of $(2k+1)$ elements. Then $<V,\cal{E}, \cal{F}>$ does not have Property S, if $\cal{E} \cap \cal{F}$ contains all the $(k+1)$ element subsets of $V$.
\end{theorem}
{\bf Proof}. Suppose $\cal{E} \cap \cal{F}$ contains all the $(k+1)$ element subsets of $V$. Any set in $\cal{E} \cap \cal{F}$ must intersect both cells of an S-partition for  $<V,\cal{E}, \cal{F}>$; however, any S-partition has one of its cells of cardinality at least $(k+1)$, and thus it must contain a set in $\cal{E} \cap \cal{F}$. \bigskip

Next we generalize to Property S, a result of Woodall \cite{Woo} for Property B .  If $\mathcal{A}$ is a family of subsets of $V$, let $\mathcal{I(A)}=\{v \in V | v \supseteq a, \hspace{.5em} \mbox{for some} \hspace{.5em} a\in \mathcal{A} \}$. Then $<V,\cal{E}, \cal{F}>$ has Property S if and only if there exists an $X$ such that $X \notin \mathcal{I(F)}$ and $V-X \notin \mathcal{I(E)}$. 
\begin{theorem} $<V,\cal{E}, \cal{F}>$ has Property S if $|\mathcal{I(E)} \cup \mathcal{I(F)}| <  2^{(|V|-1)}$.
\end{theorem}
{\bf Proof}. We claim that there is an $X$ such that $X \notin \mathcal{I(F)}$ and $V-X \notin \mathcal{I(E)}$.  If not, then for every $X \subseteq V$, either $X \in \mathcal{I(F)}$ or $V-X \in \mathcal{I(E)}$. Thus, either $X$ or $V-X$ belongs to $\mathcal{I(E)} \cup \mathcal{I(F)}$, for every $X \subseteq V$. However, this implies $|\mathcal{I(E)} \cup \mathcal{I(F)}| \geq  2^{(|V|-1)}$.

\begin{theorem}
Let V be a finite set. Then $<V,\cal{E},\cal{F}>$ has Property S if ${\sum}_{A \in \cal{E} \cup \cal{F}} {( \frac{1}{2^{|A|}})}< \frac{1}{2}$.
\end{theorem}

{\bf Proof}. Suppose ${\sum}_{A \in \mathcal{E} \cup \mathcal{F}} {(\frac{1}{2^{|A|}})< \frac{1}{2}} $. The number of subsets of $V$ which contain $A$ is $2^{(|V|-|A|)}$. Thus $|\mathcal{I(E)} \cup \mathcal{I(F)}|$$<  {\sum}_{A \in \mathcal{E} \cup \mathcal{F}}$ $2^{(|V|-|A|)}$. But ${\sum}_{A \in \mathcal{E} \cup \mathcal{F}}$ $2^{(|V|-|A|)}$$=2^{|V| } {\sum}_{A \in \mathcal{E}\cup\mathcal{F}}$$\frac{1}{2^{|A|}}<2^{(|V|-1)}$. Hence $|\mathcal{I(E)} \cup \mathcal{I(F)}| <  2^{(|V|-1)}$ and $<V,\cal{E}, \cal{F}>$ has Property S, by the previous Theorem.

If, in the last result $\mathcal{E} = \mathcal{F}$, we get the proposition in \cite{Woo}.

\begin{corollary}
If $\mathcal{F}$ is a family of subsets of a finite set $V$ and if ${\sum}_{A \in \mathcal{F} } {(\frac{1}{2^{|A|}})< \frac{1}{2}} $, then $\mathcal{F}$ has Property B.
\end{corollary}

\section{Property S and Satisfiability}

If $<V,\cal{E}, \cal{F}>$ has property S, where $\{X,V-X\}$ is an S-partition of $V$, then $X$  picks some elements from each of the subsets in $\cal{E}$ but the totality of chosen elements does not contain one of the `forbidden subsets' in $ \cal{F}$. This problem of picking elements subject to constraints is very general; we give a few examples to indicate its wide applicability. 

In propositional logic, a literal is either a statement variable or its negation, a clause is a finite collection of literals, and a conjunctive normal form (cnf) is a finite collection of clauses. The problem is to determine whether a given cnf is satisfiable; that is, does there exist an interpretation of the variables (as true or false) such that each clause contains at least one true literal.  Let $V$ be the set consisting of all the literals in the cnf, $\cal{E}$, the set of clauses of the cnf, and $ \cal{F}$, the set of the pairs consisting of the statement letters and their negations. Then the satisfiability of the cnf is easily seen to be equivalent to  $<V,\cal{E}, \cal{F}>$ has property S.

An n-coloring of a graph $G=<A,E>$ is a function $f:A\to \{1,...,n\}$, such that $f(a_1)\neq f(a_2)$, if ${a_1}E{a_2}$; $G$ is n-colorable if such a coloring exits.
Let $V$ be the pairs, $\{a,j\}$, where $a \in A$, $1\leq j\leq n$. Let $\cal{E}$ be the sets $\{\{a,1\},...,\{a,n\}\}$, where $a \in A$. Let $ \cal{F}$ consist of all the pairs, $\{\{a_1,i\},\{a_2,i\}\}$, with  ${a_1}E{a_2}$, $1\leq i \leq n$.  Then $G$ is n-colorable if and only if $<V,\cal{E}, \cal{F}>$ has property S. 

Suppose that for each vertex $a$ of a graph $G=<A,E>$, a list $L(a)$ of colors available for $a$ is given. Then a {\em list coloring} from $L$ is a proper coloring, $f$, such that $f(a) \in L(a)$; in case such a coloring exists we say that $G$ is L-list colorable. Let $\cal{E}$ be the sets $\{\{a,i\}|  i \in L(a)\}$, where $a \in A$.  Let $ \cal{F}$ consist of all the pairs, $\{\{a_1,i\},\{a_2,i\}\}$, with  ${a_1}E{a_2}$. Then $G$ is L-list colorable if and only if $<V,\cal{E}, \cal{F}>$ has property S. A graph is said to be {\em k-choosable} if it has a list coloring for every assisgnment of $k$ element lists to the vertices. (An introduction to list coloring and choosability can be found in West \cite{Wes}.)

(The ``marriage problem'') Let $S={\{S_i\}}_{i \in I}$ be an indexed family of finite sets. A system of distinct representatives (SDR) for $S$ is a one-to-one function $f:I \to \cup S$ such that $f(i) \in S_i$, $i \in I$.  Let $V$ be the set of all pairs, $\{s,i\}$, where $s \in S_i$, $i \in I$. Let $\mathcal{E}$ be all the sets $\{\{s,i\} \mid  {s \in S_i }\}$, $i \in I$. Let $ \cal{F}$ consist of the sets, $\{\{s,i\},\{s,j\}\}$, where $i \neq j$. Then $S$ has an SDR if and only if $<V,\cal{E}, \cal{F}>$ has property S.

Since so many satisfiability problems can be treated as Property S problems, it makes sense to study methods that can determine whether Property S holds.  In \cite{Cow93b}, we introduced a Tableau Method reminiscent of the Analytic Tableaux  in logic of Raymond Smullyan \cite{Smu}. The other main method, which also comes from logic, is resolution, which we turn to next.

\section{Resolution}
The Resolution proof procedure, used in logic to determine the satisfiability of conjunctive normal forms, has been generalized to provide a proof procedure for a variety of satisfiability problems by Cowen \cite{Cow91} and Kolany \cite{Kol}. 
In Kolany \cite{Kol} a notion of satisfiability on Hypergraphs was introduced which is equivalent to Property S. This enabled him to prove a far more useful resolution result than that of Cowen \cite{Cow91} whose definition of satisfiablility was too restrictive. 

\begin{definition} Let $c_1,c_2,...c_n,d,e$ be subsets of $V$. then $e$ follows from $c_1,c_2,...c_n$ by {\bf resolution} on $d$ if $d=\{v_1,v_2,...,v_n\}$, where $v_i\in c_i$, $1\le i \le n$,  and  $e=\cup_{i=1}^n c_i/v_i$; in this case, we write, $c_1,c_2,...,c_n \vdash_d e$,  
\end{definition}
If $\cal{A}$ ,$\cal{D}$ are collections of subsets of $V$, $[\cal{A}]_\mathcal{D}$ will denote the closure of $\cal{A}$ with respect to resolution on elements $d \in \mathcal{D} $. The following Theorem is elegantly proved by Kolany \cite{Kol}.

\begin{theorem} $<V,\cal{E}, \cal{F}>$ fails to have Property S if and only if $ \emptyset \in {[\cal{E}]}_{\cal{F}}$.
\end{theorem}

It is interesting to note that Woodall \cite{Woo} introduced for Property B a ``reduction process'' which is essentially the same as resolution and proved a very similar result to Kolany's Theorem (see Proposition 2 of \cite{Woo}). 

Surely $<V,\cal{E}, \cal{F}>$ has Property S if and only if $<V,\cal{F}, \cal{E}>$ has Property S. Thus, by Kolany's Theorem,  $ \emptyset \in {[\cal{E}]}_{\cal{F}}$ if and only if   $ \emptyset \in {[\cal{F}]}_{\cal{E}}$. Thus there are two distinct resolution methods which can be used in a particular case. (Of course, in the case of Property B, both methods are the same.) The next result shows how they can be combined for even greater flexibility.

\begin{theorem}  $<V,\cal{E}, \cal{F}>$ has Property S if and only if  $<V,{[\cal{E}]}_{\cal{F}}, \cal{F}>$ has Property S.
\end{theorem}
{\bf Proof}. If $<V,{[\cal{E}]}_{\cal{F}}, \cal{F}>$ has Property S, then so does $<V,\cal{E}, \cal{F}>$, since $\mathcal{E} \subset {[\mathcal{E}]}_{\mathcal{F}}$. 
Suppose that $<V,\cal{E}, \cal{F}>$ has Property S and $X$ is a transversal for $\cal{E}$ and $V-X$ is a transversal for $\cal{F}$. We claim that $X$ is also a transversal for ${[\cal{E}]}_{\mathcal{F}}$. Suppose $c_i \cap X \neq \emptyset$, $1\le i \le n$, and $d$ follows from the $c_i$ by resolution on $f$, where$f\in \mathcal{F}$, $f=\{v_1,...,v_n\}$, with $v_i \in c_i$, and $d=\cup_{i=1}^n c_i/v_i$. Then we must show that $d \cap X \neq \emptyset$, as well. 
Suppose, on the contrary that $d \cap X = \emptyset$. It follows that $c_i/v_i \cap X =\emptyset$, for $1\le i \le n$. Since $V-X$ is a transversal for $\mathcal{F}$, $v_j \in V-X$, for some j, $1\le j \le n$. Thus $v_j \not\in X$ and this, together with $c_j/v_j \cap X = \emptyset$, implies $c_j \cap X = \emptyset$; however this contradicts $c_i \cap X \neq \emptyset$, $1\le i \le n$.\bigskip

It follows that $<V,\cal{E}, \cal{F}>$ fails to have Property S if and only if $\emptyset \in {[\mathcal{F}]}_{[\mathcal{E}]_{\mathcal{F}}}$, or, if $\emptyset \in {[\mathcal{E}]}_{[\mathcal{F}]_{\mathcal{E}}}$, etc. Since this quickly becomes a typographical nightmare, we define, recursively, the following notation.
\begin{definition} $[\mathcal{F},\mathcal{E},0]=\mathcal{F}$; $[\mathcal{E},\mathcal{F},0]=\mathcal{E}$; $[\mathcal{F},\mathcal{E},n+1]={[\mathcal{F}]}_{[\mathcal{E},\mathcal{F},n]}$; $[\mathcal{E},\mathcal{F},n+1]={[\mathcal{E}]}_{[\mathcal{F},\mathcal{E},n]}$.
\end{definition}
We then have the following corollary to the previous Theorem.
\pagebreak
\begin{corollary} For $n>0$, the following statements are equivalent.
\begin{enumerate}
\item$<V,\cal{E}, \cal{F}>$ fails to have Property S 
\item  $\emptyset \in {[\mathcal{E},\mathcal{F},n]}$ 
\item $\emptyset \in {[\mathcal{F},\mathcal{E},n]}$.
\end{enumerate}
\end{corollary}
This allows resolution proofs which go back and forth; we illustrate some of these possibilities by giving three different proofs that a set $\cal{E}$ of clauses in propositional logic is unsatisfiable. We number the clauses in $\cal{E}$ for easy reference, as follows:
\begin{center} 1)$\{p,q\}$, 2)$\{p,\neg q,r\}$, 3)$\{p,\neg q,\neg r\}$, 4)$\{\neg p,q,r\}$, 5)$\{\neg p,q,\neg r\}$, 6)$\{\neg p,\neg q\}$. \end{center}
The clauses in $\cal{F}$ consist of variables and their negations:
\begin{center} A)$\{p,\neg p\}$, B)$\{q,\neg q\}$, C)$\{r,\neg r\}$. \end{center}
We annotate our proofs using the notation: $(x,y,z,.../w)$; this means that the clause was obtained from clauses labeled $x,y,z,...$ by resolving on the clause labeled $w$. Our first proof is that $\emptyset \in {[\cal{E}]}_{\cal{F}}$.
\begin{tabbing}
\hspace{1cm}\=7) $\{p,\neg q\}$  \hspace{2cm}\= (2,3/C)\\
\>8) $\{\neg p,q\}$ \>(4,5/C)\\
\>9) $\{\neg q\}$ \>  (6,7/A)\\
\>10) $\{q\}$ \> (1,8/A)\\
\>11)  $\emptyset$ \> (9,10/B) 
\end{tabbing}

Our next proof shows that $\emptyset \in {[\cal{F}]}_{\cal{E}}$.

\begin{tabbing}
\hspace{1cm}\=D) $\{p,q\}$ \hspace{2cm}\= (A,B/6)\\
\>E) $\{q,r\}$ \> (D,B,C/3)\\
\>F) $\{q\}$ \>(D,B,E/2)\\
\>G) $\{p,r\}$ \> (A,D,C/5)\\
\>H) $\{p\}$ \> (A,D,G/4)\\
\>I) $\emptyset$ \> (H,F/1)
\end{tabbing}

Our last proof demonstrates that $\emptyset \in {[\mathcal{F}]}_{[\mathcal{E}]_{\mathcal{F}}}$. Lines 12 and 13 show membership in  $[\mathcal{E}]_{\mathcal{F}}$; lines J,K show membership in ${[\mathcal{F}]}_{\mathcal{E}}$; L-N show membership in ${[\mathcal{F}]}_{[\mathcal{E}]_{\mathcal{F}}}$.
\begin{tabbing}
\hspace{1cm}\=12) $\{p,\neg q\}$ \hspace{2cm}\= (2,3/C)\\
\>13) $\{\neg p,q\}$ \> (4,5/C)\\
\> J) $\{p,q\}$ \> (A,B/6)\\
\> K) $\{\neg p, \neg q\}$ \> (A,B/1)\\
\> L) $\{p\}$ \> (A,J/13)\\
\> M) $\{\neg q\}$ \> (B,K/13)\\
\> N) $\emptyset$ \>  (L,M/12)
\end{tabbing}
\pagebreak
\bigskip
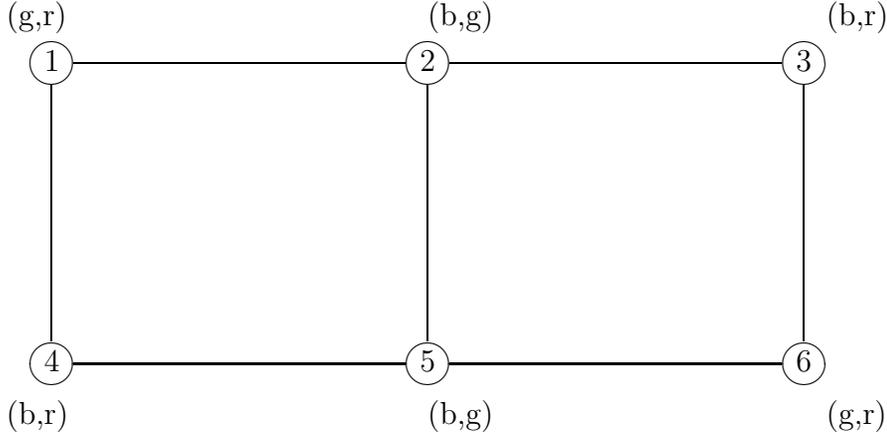
\setlength{\unitlength}{1 cm}
\begin{picture}(7,5)
\multiput(2,1)(5,0){3}{\circle{.6}}
\multiput(2,5)(5,0){3}{\circle{.6}}
\put(2.3,1){\line(1,0){4.4}}
\put(2.3,5){\line(1,0){4.4}}
\put(7.3,1){\line(1,0){4.4}}
\put(7.3,5){\line(1,0){4.4}}
\put(2,1.3){\line(0,1){3.4}}
\put(7,1.3){\line(0,1){3.4}}
\put(12,1.3){\line(0,1){3.4}}
\put(1.9,.9){4}
\put(6.9,.9){5}
\put(11.9,.9){6}
\put(1.9,4.9){1}
\put(6.9,4.9){2}
\put(11.9,4.9){3}
\put(1.4,5.5){(g,r)}
\put(7,5.5){(b,g)}
\put(12.3,5.5){(b,r)}
\put(1.4,.2){(b,r)}
\put(7,.2){(b,g)}
\put(12.3,.2){(g,r)}
\end{picture}
\bigskip

Finally we prove, by resolution, that the graph depicted above is not L-list colorable, where L is the lists shown next to the vertices. Since each list has two elements, this will show that the graph is not 2-choosable (even though it is 2-colorable). Let $\cal{E}$ be the collection of the following sets.
\begin{center} 1) $\{g_1,r_1\}$, 2) $\{b_2,g_2\}$, 3) $\{b_3,r_3\}$,  4) $\{b_4,r_4\}$, 5) $\{b_5,g_5\}$, 6) $\{g_6,r_6\}$  \end{center}
Let $\cal{F}$ be the following sets.
\begin{center} A) $\{r_1,r_4\}$, B) $\{g_1,g_2\}$, C) $\{b_4,b_5\}$, D) $\{b_2,b_5\}$, E) $\{g_2,g_5\}$, F) $\{b_2,b_3\}$,  \\ 
G) $\{r_3,r_6\}$, H) $\{g_5,g_6\}$ \end{center}

\begin{tabbing}
\hspace{1cm}\=7) $\{r_1,b_2\}$  \hspace{2cm}\= (1,2/B)\\
\>8) $\{b_3,g_6\}$ \>(3,6/G)\\
\>9) $\{b_2,b_5\}$ \>  (2,5/E)\\
\>10) $\{b_2,b_4\}$ \> (4,7/A)\\
\>11)  $\{b_3,b_5\}$ \> (5,8/H)\\
\> 12) $\{b_2\}$ \> (9,10/C) \\
\> 13) $\{b_5\}$  \> (11,12/F) \\
\> 14) $\emptyset$ \> (12,13/D)
\end{tabbing}

We leave it to the reader to provide other resolution proofs of the non list colorability.
\section{Propositional Representation}
We have seen in sections 3,4 that the satisfiability of a conjunctive normal form in propositional logic can be treated as a Property S problem. Conversely, we will show that any Property S problem can be represented as a conjunctive normal form problem. Let $<V,\cal{E}, \cal{F}>$ be a bihypergraph. For each $v \in V$, take a propositional letter $p_v$. If $E \in \cal{E}$, $E=\{e_1,...,e_n\}$, let $c_E = (p_{e_1} \vee ... \vee p_{e_n})$; if $F \in \cal{F}$, $F=\{f_1,...,f_k\}$, let $c_F=(\neg p_{f_1} \vee ... \vee \neg p_{f_k})$.
\begin{theorem}
Let $<V,\cal{E}, \cal{F}>$ be a finite bihypergraph. Then $<V,\cal{E}, \cal{F}>$ has Property S if and only if the cnf $C= \bigwedge \limits_{E \in \cal{E}} c_E \wedge \bigwedge \limits_{F \in \cal{F}} c_F$ is satisfied.
\end{theorem}
{\bf Proof}. Suppose $<V,\cal{E}, \cal{F}>$ has Property S and $\{X,V-X\}$ is an S-partition. Assign truth values to the $p_v$, $v \in V$ as follows: $p_v$ is true if and only if $v \in X$. Since $X$ is a transversal for $\cal{E}$, $X$ intersects each $E \in \cal{E}$; thus, at least one $e_i \in E$ belongs to $X$ and $p_{e_i}$ is true. Therefore $\bigwedge_{E \in \cal{E}} c_E$ is true.  If $F \in \cal{F}$, some $f_j \in F$ does not belong to $X$, since $F \not\subset X$. Therefore $p_{f_j}$ is false and $\neg p_{f_j}$ is true; hence $c_F$ is true, $F \in \cal{F}$; that is, $\bigwedge_{F \in \cal{F}} c_F$ is true. Therefore $C$ is true under the assignment.

Suppose C is satisfiable and let $I$ be a satisfying assignment for the $p_v$, $v \in V$. Let $X$ be the set of $v \in V$ such that $p_v$ is true under $I$. Suppose $E \in \cal{E}$; since $c_E$ must be true under $I$, at least one $p_{e_i}$, $e_i \in E$ must be true and hence $e_i \in X$. Thus $X$ is a transversal for $\cal{E}$. Suppose $F \in \cal{F}$; since $c_F$ is true under $I$, $\neg p_{f_j}$ must be true for at least one $f_j \in F$. Then $p_{f_j}$ is false under $I$ and so $f_j \not\in X$. Therefore $V-X$ is a transversal for $\cal{F}$. Hence $\{X,V-X\}$ is an S-partition for $<V,\cal{E}, \cal{F}>$. \bigskip

It follows that deciding whether $<V,\cal{E}, \cal{F}>$ has Property S when the sets $\cal{E}$, $\cal{F}$ consist only of pairs is polynomial-time decidable, since the conjuncts obtained all have exactly two literals and deciding CNF satisfiability in this case is known to be polynomial (see \cite{Gar}).
 
\section{Conclusion}
The various resolution techniques outlined above have yet to be tried on large examples  in logic or elsewhere. This raises several questions. For the non-logical applications such as graph coloring, how does the Property S  approach compare with other more direct techniques? In logic, how can it be determined which of the various resolution techniques outlined above for testing CNFs for satisfiability should be employed? For example, which technique should be utilized if the number of clauses is much greater than the number of variables? 

Finally, as we remarked earlier, much is known about the set-theoretic properties of Property B, in the cases of finite and infinite sets. How much of this carries over to Property S and can the results be combined with resolution to decide Property S more effectively?

\pagebreak

\end{document}